\documentclass{amsart}
\usepackage{amsmath,amsthm,amssymb,bbm}
\usepackage[bookmarksnumbered=true,pagebackref]{hyperref}

\usepackage{graphicx}

\newtheorem{thm}{Theorem}[section]
\newtheorem{lem}[thm]{Lemma}

\theoremstyle{definition}

\theoremstyle{remark}
\newtheorem{rem}[thm]{Remark}

\numberwithin{equation}{section}

\newcommand{\D}[1]{\mathop{\mathrm{d}#1}}
\DeclareMathOperator\supp{supp}

\begin{document}

\title[Spectral upper bound for the torsion function]{Spectral upper bound for the torsion function of symmetric stable processes}


\author{Hugo Panzo}
\address{Technion -- Israel Institute of Technology\\
Haifa 32000, Israel}
\email{panzo@campus.technion.ac.il}
\thanks{Supported at the Technion by a Zuckerman Fellowship.}



\date{\today}

\dedicatory{}


\begin{abstract}
We prove a spectral upper bound for the torsion function of symmetric stable processes that holds for convex domains in $\mathbb{R}^d$. Our bound is explicit and captures the correct order of growth in $d$, improving upon the existing results of Giorgi and Smits \cite{Smits} and Biswas and L\H{o}rinczi \cite{Biswas}. Along the way, we make progress towards a torsion analogue of Chen and Song's \cite{two_sided} two-sided eigenvalue estimates for subordinate Brownian motion. 
\end{abstract}

\maketitle

\section{Introduction}\label{sec:intro}

Suppose $D\subset\mathbb{R}^d$ is a domain (nonempty connected open set) and consider the $d$-dimensional Brownian motion $W=(W_t:t\geq 0)$ starting at $x\in D$ and running at twice the usual speed until the first exit time $\tau_D:=\inf\{t\geq 0:W_t\notin D\}$ where we take $\inf\{\emptyset\}=\infty$. It is well known that for any starting point $x\in D$, the right tail of the exit time has an exponential rate of decay given by the \emph{principal eigenvalue} of the Laplacian $\Delta$ on $D$ with Dirichlet boundary conditions. More precisely, it follows from \cite[Theorem 3.1.2]{Sznitman} that for all $x\in D$ we have
\begin{equation}\label{eq:principal}
-\lim_{t\to\infty}\frac{1}{t}\log\mathbb{P}_x\left(\tau_D >t\right)=\lambda_D
\end{equation}
where
\begin{equation}\label{eq:Rayleigh_quotient}
\lambda_D:=\inf_{\substack{\varphi\in C_c^\infty(D)\\ \varphi\not\equiv 0}}\frac{\int_D|\nabla \varphi|^2\D{x}}{\int_D \varphi^2\D{x}}.
\end{equation}
Here $C_c^\infty(D)$ denotes the space of smooth functions on $D$ with compact support. In general, $\lambda_D$ needn't be an actual eigenvalue but merely the bottom of the spectrum.

The expected exit time from $D$ as a function of the starting position, that is 
\begin{equation}\label{eq:torsion}
u_D(x):=\mathbb{E}_x\left[\tau_D\right],~x\in D,
\end{equation}
is called the \emph{torsion function} of the domain $D$. A result of Burkholder \cite[ Equation 3.13]{Burkholder} shows that if $u_D(x)<\infty$ for some $x\in D$, then $u_D(x)<\infty$ for all $x\in D$. This in turn implies that $u_D$ is the unique solution of the boundary value problem
\begin{equation}\label{eq:PDE}
\left\{
\begin{aligned}
-\Delta u&=1 & &\text{in }D\\ 
u&=0 & &\text{on }\partial D
\end{aligned}
\right.
\end{equation}
where the Dirichlet boundary condition is understood in the Sobolev sense; see \cite{vdB_Carroll,vdB_1}. The integral of $u_D$ over $D$, namely $\|u_D\|_1$, is known as the \emph{torsional rigidity} of $D$ and it can be used to quantify the resistance to twisting of a beam with cross section $D$. For other applications of the torsion function, the reader is directed to the solution of the classical de Saint-Venant problem \cite{Polya} and more recent results related to Anderson localization \cite{landscape} and Hermite-Hadamard inequalities \cite{Hermite_Hadamard}.

The scaling property of Brownian motion implies that if we scale $D$ by a factor of $s>0$, then the new principal eigenvalue and torsion function are $\frac{1}{s^2}\lambda_D$ and $s^2\,u_D(\frac{x}{s})$, respectively. Since $D$ contains a small enough ball, we know that $\lambda_D<\infty$ and $\|u_D\|_\infty>0$. If we also assume that either $\lambda_D>0$ or $\|u_D\|_\infty<\infty$ holds, then the product $\lambda_D\,\|u_D\|_\infty$ is well-defined and scale invariant. In fact, a result of van den Berg and Carroll \cite{vdB_Carroll} shows that $\lambda_D>0$ if and only if $\|u_D\|_\infty<\infty$. Changing the shape of $D$, however, can have unequal competing effects on $\lambda_D$ and $\|u_D\|_\infty$ which make $\lambda_D\,\|u_D\|_\infty$ an interesting quantity to study. See \cite{vdB_capacity,vdB_Santalo} for recent work on other ``shape functionals'' that involve the principal eigenvalue or torsional rigidity.

Finding upper and lower bounds on the product $\lambda_D\,\|u_D\|_\infty$ which hold for various classes of domains $D$ has been a topic of active research.  Define the constants
\begin{equation*}
\begin{split}
m_d&=\inf\{\lambda_D\,\|u_D\|_\infty:D\subset\mathbb{R}^d\text{ is a domain with }\lambda_D>0\},\\
M_d&=\sup\{\lambda_D\,\|u_D\|_\infty:D\subset\mathbb{R}^d\text{ is a domain with }\lambda_D>0\}.
\end{split}
\end{equation*}
Then for any domain $D\subset\mathbb{R}^d$ with $\lambda_D>0$ we have by definition  
\begin{equation}\label{eq:torsion_bound}
m_d\leq\lambda_D\,\|u_D\|_\infty\leq M_d.
\end{equation}
We can also get a version of \eqref{eq:torsion_bound} which holds for any domain regardless of whether or not $\lambda_D>0$. There are several ways to show that $\lambda_D\,\|u_D\|_\infty\geq 1$ for bounded domains; see \cite{Rodrigo,vdB_Carroll,Henrot}. Extending this inequality to general domains with $\lambda_D>0$ is also a straightforward matter, hence $M_d\geq m_d\geq 1$ for all $d$. Now the fact mentioned above that $\lambda_D=0$ if and only if $\|u_D\|_\infty=\infty$ implies that the inequality
\begin{equation}\label{eq:spectral_bound}
\frac{m_d}{\lambda_D}\leq\|u_D\|_\infty\leq\frac{M_d}{\lambda_D}
\end{equation}
holds for any domain $D\subset\mathbb{R}^d$ where in case $\lambda_D=0$ we take $1/0=\infty$. 

Inequalities such as \eqref{eq:torsion_bound} and \eqref{eq:spectral_bound} are known as \emph{spectral bounds for the torsion function} and have appeared in \cite{Rodrigo,vdB_Carroll,Smits,vdB_1,Henrot,vdB_2,Vogt,Biswas,improved_Vogt}, sometimes with different constants that reflect various further assumptions on $D$. For instance, there has also been interest in the analogous inequalities for convex domains which feature the similarly defined constants $c_d$ and $C_d$ given by
\begin{equation}\label{eq:C_constant}
\begin{split}
c_d&=\inf\{\lambda_D\,\|u_D\|_\infty:D\subset\mathbb{R}^d\text{ is a convex domain with }\lambda_D>0\},\\
C_d&=\sup\{\lambda_D\,\|u_D\|_\infty:D\subset\mathbb{R}^d\text{ is a convex domain with }\lambda_D>0\}.
\end{split}
\end{equation}
See \cite{Biswas,improved_Vogt} for a generalization that involves the $p$-th moment of the exit time. A recent application of the upper bound in \eqref{eq:spectral_bound} to extremal problems related to the conformal Skorokhod embedding can be found in \cite{CSEP}.

Clearly $m_d\leq c_d\leq C_d\leq M_d$, and since domains in $\mathbb{R}$ are just intervals, straightforward calculations show that $m_1=c_1=C_1=M_1=\frac{\pi^2}{8}$. For the nontrivial case $d\geq 2$, the papers \cite{Henrot,vdB_2} construct a family of ``Swiss cheese'' domains to show that $m_d=1$, while it follows from a result of Payne \cite{Payne} that $c_d=\frac{\pi^2}{8}$ and this is attained by an infinite slab in $\mathbb{R}^d$. On the other hand, much less is known about $C_d$ and $M_d$ for $d\geq 2$. However, a recent result of H. Vogt \cite{Vogt} shows that 
\begin{equation}\label{eq:Vogt_C}
M_d\leq\underbrace{\frac{d}{8}+\frac{\sqrt{d}}{4}\sqrt{5\left(1+\frac{1}{4}\log 2\right)}+1}_{\displaystyle \mathcal{V}_d}
\end{equation}
and that $\frac{d}{8}$ is the correct leading term for large $d$. Vogt's bound $\mathcal{V}_d$ is the latest in a series of explicit upper bounds which have been improved upon over the years; see \cite{Rodrigo,vdB_Carroll,Smits} for some earlier variants. See also \cite{improved_Vogt} for a numerical improvement of Vogt's bound. We refer to \cite{Henrot,improved_Vogt} for some theorems and conjectures regarding sharp upper bounds and the existence of extremal domains which attain them. 

\subsection{\texorpdfstring{Other conventions for the product $\lambda_D\,\|u_D\|_\infty$}{Other conventions for the product}}

From the scaling considerations discussed above, it follows that the product $\lambda_D\,\|u_D\|_\infty$ is the same as that for Brownian motion running at its usual speed, or, for that matter, at any constant multiple of its usual speed. However, some authors \cite{Rodrigo,improved_Vogt} have considered Brownian motion running at twice its usual speed in \eqref{eq:principal} and at its usual speed in \eqref{eq:torsion} and this results in the product $\lambda_D\,\|u_D\|_\infty$ being twice that of what one obtains from using either speed consistently in both \eqref{eq:principal} and \eqref{eq:torsion}. This factor of $2$ also shows up in their versions of the constants $m_d$, $c_d$, $C_d$, and $M_d$ as well as the bound $\mathcal{V}_d$. In this paper we always consider the product $\lambda_D\,\|u_D\|_\infty$ obtained by using the process run at the same speed in both \eqref{eq:principal} and \eqref{eq:torsion} so there is never any ambiguity. However, the reader is advised to exercise caution when consulting the literature as there is no consensus on which product to consider. 

\section{Main results}

The goal of this paper is to improve the existing spectral upper bound for the torsion function of the \emph{symmetric stable processes}. By symmetric we mean rotationally symmetric or isotropic. Recall that for $0<\alpha\leq 2$, the $d$-dimensional symmetric $\alpha$-stable process is the L\'{e}vy process $X=(X_t:t\geq 0)$ with characteristic function 
\begin{equation}\label{eq:stable}
\mathbb{E}\left[e^{i\xi\cdot X_1}\right]=e^{-|\xi|^\alpha},~\xi\in\mathbb{R}^d.
\end{equation}
Note that when $\alpha=2$, the process $X$ is simply Brownian motion run at twice the usual speed, namely, our process $W$.

The probabilistic definitions of $\lambda_D$ and $u_D$ appearing in \eqref{eq:principal} and \eqref{eq:torsion} are perfectly valid when $W$ is replaced by $X$. Moreover, the scale invariance of the product $\lambda_D\,\|u_D\|_\infty$ follows from the $\alpha$-self-similarity of the paths of $X$, just like in the Brownian case. However, the Laplacian and Dirichlet form appearing in the PDE \eqref{eq:PDE} and variational \eqref{eq:Rayleigh_quotient} formulations of $u_D$ and $\lambda_D$ need to be replaced by the Dirichlet fractional Laplace operator and its associated quadratic form; see \cite{Frank,fractional_survey} and references therein. 

To avoid ambiguity, from now on we attach superscripts to $\lambda_D$ and $u_D$ which indicate the corresponding process. In \cite[Theorem 3.1]{Smits}, Giorgi and Smits prove that there exist constants $M(d,\alpha)$ such that 
\begin{equation}\label{eq:Smits}
1\leq\lambda_D^X\,\left\|u_D^X\right\|_\infty\leq M(d,\alpha)
\end{equation}
holds for any domain $D\subset\mathbb{R}^d$ where the transition density of $X$ killed upon exiting $D$ admits a Hilbert-Schmidt expansion. For instance, this holds when $D$ has finite volume. They mention in Remark 3.3 that there is room for improvement in the upper bound since their proof can only provide constants $M(d,\alpha)$ that grow superexponentially in $d$. More recently, similar nonexplicit bounds have been established by Biswas and L\H{o}rinczi \cite[Corollary 3.3]{Biswas} which hold for a more general class of L\'{e}vy processes under the additional constraint that $D$ is a bounded convex domain.

The first main result of this paper is the following improvement of \eqref{eq:Smits}, that when paired with \eqref{eq:Vogt_C}, leads to an explicit spectral upper bound which holds for all convex domains $D\subset\mathbb{R}^d$ and captures the correct order of growth in $d$.
\begin{thm}\label{thm:main}
Let $X$ be a $d$-dimensional symmetric $\alpha$-stable process with $0<\alpha\leq 2$ and let $C_d$ be defined by \eqref{eq:C_constant}. If $D\subset\mathbb{R}^d$ is a convex domain then 
\begin{equation}\label{eq:main}
\frac{1}{\lambda_D^X}\leq\left\|u_D^X\right\|_\infty\leq\frac{4}{\alpha\,\Gamma(\alpha/2)}\frac{C_d^{\alpha/2}}{\lambda_D^X}
\end{equation}
where we take $1/0=\infty$. In particular, Vogt's result \eqref{eq:Vogt_C} gives the explicit bound
\begin{equation*}
\left\|u_D^X\right\|_\infty\leq\frac{4}{\alpha\,\Gamma(\alpha/2)}\frac{\mathcal{V}_d^{\alpha/2}}{\lambda_D^X}.
\end{equation*}
\end{thm}
\begin{rem}
As discussed in the paragraph preceding \eqref{eq:Vogt_C}, the $1$ appearing in the lower bound of \eqref{eq:main} can be sharpened to $\frac{\pi^2}{8}$ when $\alpha=2$ and this is attained by an infinite slab in $\mathbb{R}^d$. An interesting question is whether the analogous statement also holds when $0<\alpha<2$, that is, can the $1$ be sharpened to $\lambda_I^X\,\|u_I^X\|_\infty$ where $I$ is the interval $(-1,1)$? Note that while $\|u_I^X\|_\infty$ is known explicitly \eqref{eq:ball_sup}, only (rather precise) numerical estimates are available for $\lambda_I^X$; see \cite{ball_lambda}.
\end{rem}
\begin{rem}
When $\lambda_D^X>0$, putting $\alpha=2$ in \eqref{eq:main} results in twice the upper bound we get from the definition of $C_d$ alone. This is an artefact of the repeated resurrection technique used in the proof of Lemma \ref{lem:resurrect}. It's not clear if a refinement of this technique might lead to an upper bound which becomes sharp as $\alpha\nearrow 2$.
\end{rem}

For two functions $f,g:\mathbb{N}\to(0,\infty)$, we write $f(d)\asymp g(d)$ as $d\to\infty$ to indicate that both $f(d)=O\big(g(d)\big)$ and $g(d)=O\big(f(d)\big)$ as $d\to\infty$. The following lemma, whose proof is postponed until Section \ref{sec:proofs}, is the symmetric stable process counterpart of Lemma 1.4 in \cite{Vogt}. As alluded to in \cite[Remark 3.3]{Smits}, the lemma confirms that for the $d$-dimensional unit ball $B$, the product $\lambda_B^X\,\|u_B^X\|_\infty\asymp d^{\alpha/2}$ as $d\to\infty$. Since we also have $\mathcal{V}_d^{\alpha/2}\asymp d^{\alpha/2}$ as $d\to\infty$, the upper bound in Theorem \ref{thm:main} captures the correct order of growth in $d$.
\begin{lem}\label{lem:comparison}
Let $X$ be a $d$-dimensional symmetric $\alpha$-stable process with $0<\alpha<2$ and let $B$ be the $d$-dimensional unit ball. Then for each $\alpha$, there exist constants $C_1(\alpha)\in(0,\infty)$ and $C_2(\alpha)\in (0,\infty)$ such that 
\[
\frac{1}{2}\frac{2^{-3\alpha/2}}{\Gamma(1+\alpha/2)}d^{\alpha/2}\leq\lambda_B^X\, \left\|u_B^X\right\|_\infty\leq\frac{2^{-3\alpha/2}}{\Gamma(1+\alpha/2)}d^{\alpha/2}+C_1(\alpha)\, d^{\alpha/6}+C_2(\alpha)\, d^{\alpha-1}
\]
holds for all $d=1,2,\dots$.
\end{lem}

When first attempting to prove a result like Theorem \ref{thm:main}, one is inclined to try and adapt Vogt's method \cite{Vogt} from the Brownian case to that of stable processes. This approach soon becomes problematic since his argument uses the fact that the Dirichlet heat semigroup is dominated by the free heat semigroup on $\mathbb{R}^d$ in an essential way. More specifically, the method of weighted estimates that he employs has no hope of working with the free stable semigroup (at least not with exponential weights) due to the heavy-tailed nature of $\alpha$-stable processes with $0<\alpha<2$.

We take an alternative approach which uses Vogt's result ``off the shelf" by first realizing the symmetric $\alpha$-stable process $X$ as a subordinate Brownian motion and then using the potential theory available for such processes to carry over the upper bound \eqref{eq:Vogt_C} from the Brownian case. This still leaves open the question of whether a more direct method can be found to produce an upper bound that captures the correct order of growth in the dimension $d$, or better still, an upper bound that is asymptotically sharp like that of Vogt in the $\alpha=2$ case. 

The proof of \eqref{eq:main} in Theorem \ref{thm:main} has two main ingredients. In order to describe these components, we first recall some basic facts about subordinators and subordinate Brownian motion. The reader can consult \cite{subordinators,Bogdan_book} for more details. A (possibly killed) \emph{subordinator} is an increasing L\'{e}vy process $S=(S_t:t\geq 0)$ taking values in $[0,\infty]$ with $S_0=0$ and with $\infty$ serving as the cemetery state, if any. The \emph{Laplace exponent} $\Phi$ of $S$ is defined by 
\begin{equation*}
\Phi(\lambda)=-\log\mathbb{E}\big[\exp\left(-\lambda S_1\right)\big],~\lambda\geq 0
\end{equation*}
where by convention $e^{-\lambda\infty}=0$ for all $\lambda\geq 0$. Moreover, $\Phi$ can be written as
\begin{equation}\label{eq:LK_rep}
\Phi(\lambda)=\mathbbm{k}+\mathbbm{d}\,\lambda+\int_{(0,\infty)}\left(1-e^{-\lambda t}\right)\Pi(\D{t}),~\lambda\geq 0
\end{equation}
for unique $\mathbbm{k},\mathbbm{d}\geq 0$ and measure $\Pi$ on $(0,\infty)$ satisfying 
\[
\int_{(0,\infty)}(1\wedge t)\Pi(\D{t})<\infty.
\]

The constants $\mathbbm{k}$ and $\mathbbm{d}$ are called the \emph{killing rate} and \emph{drift}, respectively, and $\Pi$ the \emph{L\'{e}vy measure} of the subordinator $S$. We say that a subordinator is \emph{unkilled} if $\mathbbm{k}=0$. A subordinator with drift or infinite L\'{e}vy measure, that is, with $\mathbbm{d}>0$ or $\Pi\big((0,\infty)\big)=\infty$, has paths that are almost surely strictly increasing. It follows from the L\'{e}vy--Khintchine representation \eqref{eq:LK_rep} that $\Phi'$ is a completely monotone function. In particular, $\Phi$ is increasing and continuous on $[0,\infty)$. A subordinator is called a \emph{special subordinator} if its \emph{conjugate Laplace exponent} $\Phi^*(\lambda):=\lambda/\Phi(\lambda)$ is also the Laplace exponent of a subordinator. 

Additionally, the \emph{potential measure} $V$ of a subordinator $S$ is defined by
\begin{equation*}
V(A)=\mathbb{E}\left[\int_0^\infty \mathbbm{1}_{\{S_t\in A\}}\D{t}\right]
\end{equation*}
where $A\subset[0,\infty)$ is a Borel set. In other words, $V(A)$ is the expected time that $S$ spends in the set $A$. Furthermore, the Laplace transform of $V$ is related to $\Phi$ via
\begin{equation}\label{eq:Lap_link}
\int_{[0,\infty)} e^{-\lambda x}V(\D{x})=\frac{1}{\Phi(\lambda)},~\lambda>0.
\end{equation}
With a slight abuse of notation, we also use $V$ to denote the distribution function of the potential measure $V$, namely, $V(x):=V\big([0,x]\big)$. In this case we refer to $V$ as the \emph{renewal function} of $S$. Note that $\mathbbm{k}=0$ implies that $V$ is unbounded.

Suppose $S$ is an unkilled subordinator with Laplace exponent $\Phi$ and let $W$ be an independent Brownian motion in $\mathbb{R}^d$ running at twice the usual speed. We use $Y=(Y_t:t\geq 0)$ to denote the \emph{subordinate Brownian motion} which is defined by $Y_t=W_{S_t}$. It is well known that $Y$ is a L\'{e}vy process and has characteristic function
\begin{equation}\label{eq:subordinate}
\mathbb{E}\left[e^{i\xi\cdot Y_1}\right]=e^{-\Phi\left(|\xi|^2\right)},~\xi\in\mathbb{R}^d.
\end{equation}

The first component used in the proof of Theorem \ref{thm:main} is the upper bound from Chen and Song's \cite{two_sided} two-sided eigenvalue estimates for subordinate processes. Under certain conditions, these estimates allow us to bound $\lambda_D^Y$ using $\Phi(\lambda_D^W)$. More specifically, their Theorem 4.5, Remark 3.5, and the note added in proof imply that
\begin{equation}\label{eq:two_sided} 
\frac{1}{2}\Phi\left(\lambda_D^W\right)\leq \lambda_D^Y \leq \Phi\left(\lambda_D^W\right)
\end{equation}
for any bounded convex $D\subset\mathbb{R}^d$ such that the transition density of $Y$ killed upon exiting $D$ admits a Hilbert-Schmidt expansion. In particular, this inequality applies when $Y$ is a symmetric $\alpha$-stable process. 

The second component, which is of interest in its own right, is an upper bound on $\|u_D^Y\|_\infty$ that appears in the following theorem which we prove in Section \ref{sec:Theorem_2}. This theorem can be seen as one side of a torsion analogue of the Chen and Song two-sided estimate \eqref{eq:two_sided}; see Remark \ref{rem:torsion_analogue}.

\begin{thm}\label{thm:torsion_bounds}
Suppose $S$ is an unkilled special subordinator with drift or infinite L\'{e}vy measure and let $\Phi$ be its Laplace exponent and $V$ be the distribution function of its potential measure. If $D\subset\mathbb{R}^d$ is a convex domain then we have
\begin{equation}\label{eq:CS_torsion}
\frac{1}{\Phi\left(\lambda_D^W\right)}\leq\left\|u_D^Y\right\|_\infty\leq 2\,V\Big(\left\|u_D^W\right\|_\infty\Big)
\end{equation}
where $1/\Phi(0)=\displaystyle\lim_{\lambda\searrow 0}1/\Phi(\lambda)=\infty$ and $V(\infty)=\displaystyle\lim_{x\nearrow\infty}V(x)=\infty$. In particular, if either $\lambda_D^W=0$ or $\left\|u_D^W\right\|_\infty=\infty$ hold, then \eqref{eq:CS_torsion} degenerates to $\infty\leq\left\|u_D^Y\right\|_\infty\leq\infty$.
\end{thm}
\begin{rem}\label{rem:torsion_analogue}
A natural question is whether the torsion analogue of \eqref{eq:two_sided} holds, that is, can we can take $V(\|u_D^W\|_\infty)$ as a lower bound for $\|u_D^Y\|_\infty$? Evidence that this may be true, perhaps with further assumptions on $D$ or $S$, can be seen from a straightforward calculation using \eqref{eq:renewal_function} and \eqref{eq:tor_bound} which shows that this lower bound holds for all dimensions when $D$ is a ball and $Y$ is any symmetric $\alpha$-stable process.
\end{rem}

\subsection{More general domains}

The convexity assumption in our theorems can be relaxed to a uniform Lipschitz condition at the expense of introducing a nonexplicit constant into the upper bounds. More specifically, the Lipschitz condition would satisfy the hypotheses of Lemma \ref{lem:Jensen} while also implying an exterior cone condition under which Lemma \ref{lem:resurrect} would hold with the $2$ replaced by a nonexplicit constant; see \cite[Section 2]{Lip_cone} and the paragraph after \eqref{eq:geo}. The upper bound of \eqref{eq:two_sided} would also continue to hold under this exterior cone condition. However, since the explicit nature of the upper bounds in Theorems \ref{thm:main} and \ref{thm:torsion_bounds} is one of the main novelties of our results, we state and prove these theorems with a convexity assumption.

\section{Killed subordinate and subordinate killed Brownian motion}\label{sec:Theorem_2}

In order to prove Theorem \ref{thm:torsion_bounds}, we need to recall some further results about special subordinators and subordinate Brownian motion. First of all, a necessary and sufficient condition for a subordinator $S$ to be a special subordinator is that its potential measure $V$ can be written as
\begin{equation}\label{eq:special}
V(\D{t})=c\,\delta_0(\D{t})+v(t)\D{t},~t\geq 0
\end{equation}
for some $c\geq 0$ and some decreasing function $v:(0,\infty)\to(0,\infty)$ satisfying $\int_0^1 v(t)\D{t}<\infty$; see \cite[Theorem 5.1]{Bogdan_book}. Note that if \eqref{eq:special} holds, then the renewal function of $S$ is continuous and concave. Furthermore, if we insist that $S$ has drift or infinite L\'{e}vy measure, then the resulting strict monotonicity of the paths of $S$ implies that $V$ is atomless, hence $c=0$ in this case. We also note that the range of a subordinator, when restricted to $[0,\infty)$, is bounded almost surely if $\mathbbm{k}>0$. Otherwise, it is unbounded almost surely. In particular, the first passage time of an unkilled subordinator over any level is finite almost surely.

From now on, we assume that $S$ is an unkilled subordinator, that $W$ is an independent Brownian motion in $\mathbb{R}^d$ running at twice the usual speed, and that $Y$ is $W$ subordinated by $S$. Let $D\subset\mathbb{R}^d$ be a domain and define the \emph{killed Brownian motion} $W^D=\left(W^D_t:t\geq 0\right)$ by $W^D_t=W_t$ for $t<\tau_D^W$ and $W^D_t=\partial$ for $t\geq\tau_D^W$ with $\partial$ denoting the cemetery state. The \emph{killed subordinate Brownian motion} $Y^D=\left(Y^D_t:t\geq 0\right)$ is defined analogously by $Y^D_t=Y_t$ for $t<\tau_D^Y$ and $Y^D_t=\partial$ for $t\geq\tau_D^Y$ where $\tau_D^Y:=\inf\{t\geq 0:Y_t\notin D\}$. Notice that $Y^D$ results from first subordinating and then killing $W$. Switching the order of this procedure results in the \emph{subordinate killed Brownian motion} $Z^D=\left(Z^D_t:t\geq 0\right)$ which is defined by $Z^D_t=W^D_{S_t}$. Defining $\tau_D^Z$ similarly to $\tau_D^W$ and $\tau_D^Y$, we see that $\tau_D^Z=\inf\{t\geq 0:S_t\geq \tau_D^W\}$. The respective generators of $W^D$, $Y^D$, and $Z^D$ are the \emph{Dirichlet Laplacian}, the \emph{Dirichlet fractional Laplacian}, and the \emph{fractional Dirichlet Laplacian}; see \cite[Section 8]{fractional_survey}. 

The basic idea behind the proof of Theorem \ref{thm:torsion_bounds} is to exploit the close relationship between these three processes in order to compare both $\tau_D^W$ and $\tau_D^Y$ with $\tau_D^Z$, and consequently, with each other. This is a common technique in the potential theory of subordinate processes; see \cite{Bogdan_book} and references therein. 

As a first step in this scheme, we compare the expectations of $\tau_D^W$ and $\tau_D^Z$ in the following lemma.
\begin{lem}\label{lem:Jensen}
Suppose $S$ is an unkilled special subordinator with drift or infinite L\'{e}vy measure and let $\Phi$ be its Laplace exponent and $V$ be the distribution function of its potential measure. If $D\subset\mathbb{R}^d$ is a bounded Lipschitz domain then
\[
\frac{1}{\Phi\left(\lambda_D^W\right)}\leq\sup_{x\in D}\mathbb{E}_x\left[\tau_D^Z\right]\leq V\Big(\left\|u_D^W\right\|_\infty\Big).
\]
\end{lem}
\begin{proof}
The fact that $S$ is a special subordinator means that its potential measure has representation \eqref{eq:special}, hence $V$ is concave. Moreover, $c=0$ since $S$ has drift or infinite L\'{e}vy measure. Since $D$ is a bounded Lipschitz domain, the semigroup of $W^D$ is intrinsically ultracontractive. This implies that the Green function of killed Brownian motion subordinated by the special subordinator $S$ can be written as 
\begin{equation*}
G_D^Z(x,y)=\int_0^\infty p_D^W(t,x,y)v(t)\D{t},~x,y\in D
\end{equation*}
where $p_D^W$ is the transition density of the killed Brownian motion; see \cite[Section 5.5]{Bogdan_book}. Hence for $x\in D$, we can use Tonelli's theorem to write 
\begin{equation}\label{eq:transition}
\mathbb{E}_x\left[\tau_D^Z\right]=\int_D G_D^Z(x,y)\D{y}=\int_0^\infty\int_D p_D^W(t,x,y)\D{y}v(t)\D{t},
\end{equation}
and using integration by parts on the right-hand side of \eqref{eq:transition} yields
\begin{equation}\label{eq:Z_identity}
\mathbb{E}_x\left[\tau_D^Z\right]=\int_0^\infty \mathbb{P}_x\left(\tau_D^W>t\right)v(t)\D{t}=\mathbb{E}_x\left[V\left(\tau_D^W\right)\right].
\end{equation}
Now an application of Jensen's inequality to \eqref{eq:Z_identity} allows us to write
\begin{equation*}
\mathbb{E}_x\left[V\left(\tau_D^W\right)\right]=\mathbb{E}_x\left[\tau_D^Z\right]\leq V\left(\mathbb{E}_x\left[\tau_D^W\right]\right),
\end{equation*}
and noting that $V$ is increasing leads to
\begin{equation}\label{eq:Z_upper}
\sup_{x\in D}\mathbb{E}_x\left[V\left(\tau_D^W\right)\right]=\sup_{x\in D}\mathbb{E}_x\left[\tau_D^Z\right]\leq V\Big(\left\|u_D^W\right\|_\infty\Big).
\end{equation}
This proves the upper bound of Lemma \ref{lem:Jensen}.

Next we prove the lower bound. The boundedness of $D$ implies that for all $t>0$, the killed transition density $p_D^W$ has the Hilbert-Schmidt expansion 
\begin{equation*}
p_D^W(t,x,y)=\sum_{n=1}^\infty e^{-\lambda_n t}\phi_n(x)\phi_n(y),~x,y\in D
\end{equation*}
where the eigenfunctions $\{\phi_n\}_{n\geq 1}$ form an orthonormal basis for $L^2(D)$ with corresponding eigenvalues $\{\lambda_n\}_{n\geq 1}$. In particular, $\lambda_1=\lambda_D^W$ and $\phi_1$ can be taken positive so that $\phi_1/\|\phi_1\|_1$ is a probability density on $D$; see \cite[Section 4.1]{Bogdan_book}. Moreover, the series converges absolutely for all $t>0$, as seen from the Cauchy--Schwarz bound  
\begin{equation}\label{eq:Cauchy_Schwarz}
\sum_{n=1}^\infty e^{-\lambda_n t}\big|\phi_n(x)\phi_n(y)\big|\leq\sqrt{p_D^W(t,x,x)\,p_D^W(t,y,y)}\leq p^W(t,0,0)
\end{equation}
where $p^W$ is the free heat kernel on $\mathbb{R}^d$.

Note that taking a supremum over $D$ dominates averaging over $D$ with the probability density $\phi_1/\|\phi_1\|_1$. Hence we can use the equality of the right-hand sides of \eqref{eq:Z_identity} and \eqref{eq:transition}, then justify an application of Fubini's theorem with \eqref{eq:Cauchy_Schwarz}, and finally use orthogonality and \eqref{eq:Lap_link} to write
\begin{align*}
\sup_{x\in D}\mathbb{E}_x\left[V\left(\tau_D^W\right)\right]&\geq \int_D\mathbb{E}_x\left[V\left(\tau_D^W\right)\right]\frac{\phi_1(x)}{\|\phi_1\|_1}\D{x}\\
&=\int_0^\infty\int_D\int_D \sum_{n=1}^\infty e^{-\lambda_n t}\phi_n(x)\phi_n(y)\,\frac{\phi_1(x)}{\|\phi_1\|_1}\D{x}\D{y}v(t)\D{t}\\
&=\int_0^\infty e^{-\lambda_1 t}v(t)\D{t}\\
&=\frac{1}{\Phi\left(\lambda_D^W\right)}.
\end{align*}
Combining this inequality with \eqref{eq:Z_upper} proves the lower bound of Lemma \ref{lem:Jensen}.
\end{proof}

\subsection{\texorpdfstring{Constructing $\tau_D^Y$ via repeated resurrections of $Z^D$}{Constructing the Y exit time via resurrections of Z}}\label{sec:resurrect}

The second step is to compare $\tau_D^Y$ with $\tau_D^Z$. We do this by constructing $\tau_D^Y$ through repeated resurrections of $Z^D$; see \cite{killed_subordinate,Hurd} for other implementations of this procedure. To illustrate this idea, we start with the observation that $\tau_D^Y=\tau_D^Z$ on the event $\{W_{S_{\tau_D^Z}}\notin D\}$. This event occurs with probability $1$ if $S$ passes over $\tau_D^W$ continuously, for in that case $S_{\tau_D^Z}=\tau_D^W$. However, if $S$ passes over $\tau_D^W$ via a jump, then there is positive probability that $W$ may have wandered back into $D$ during the overshoot which would result in $W_{S_{\tau_D^Z}}\in D$. In this case we apply the appropriate Markov shift operators to the paths of $W$ and $S$ and then restart $Z^D$. This procedure is repeated while keeping track of each resulting $\tau_D^Z$ until finally $W_{S_{\tau_D^Z}}\notin D$. Now $\tau_D^Y$ can be obtained by summing these $\tau_D^Z$. 

The following argument uses this idea in a more precise way to establish an upper bound on the expectation of $\tau_D^Y$. First we need to define two increasing sequences of stopping times $\{\tau_n\}_{n\geq 0}$ and $\{\sigma_n\}_{n\geq 0}$ which correspond to the times at which the resurrected versions of the processes $W^D$ and $Z^D$ exit $D$. Applying the strong Markov property to $Y$ at the times $\{\sigma_n\}_{n\geq 0}$ is an essential step in our argument and this can be justified through the construction of an auxiliary filtration that contains the natural filtration of $Y$ and with respect to which $Y$ is a strong Markov process and $\{\sigma_n\}_{n\geq 0}$ are stopping times, refer to \cite[Section 2]{killed_subordinate} for more details.

\begin{figure}
\includegraphics[scale=.6]{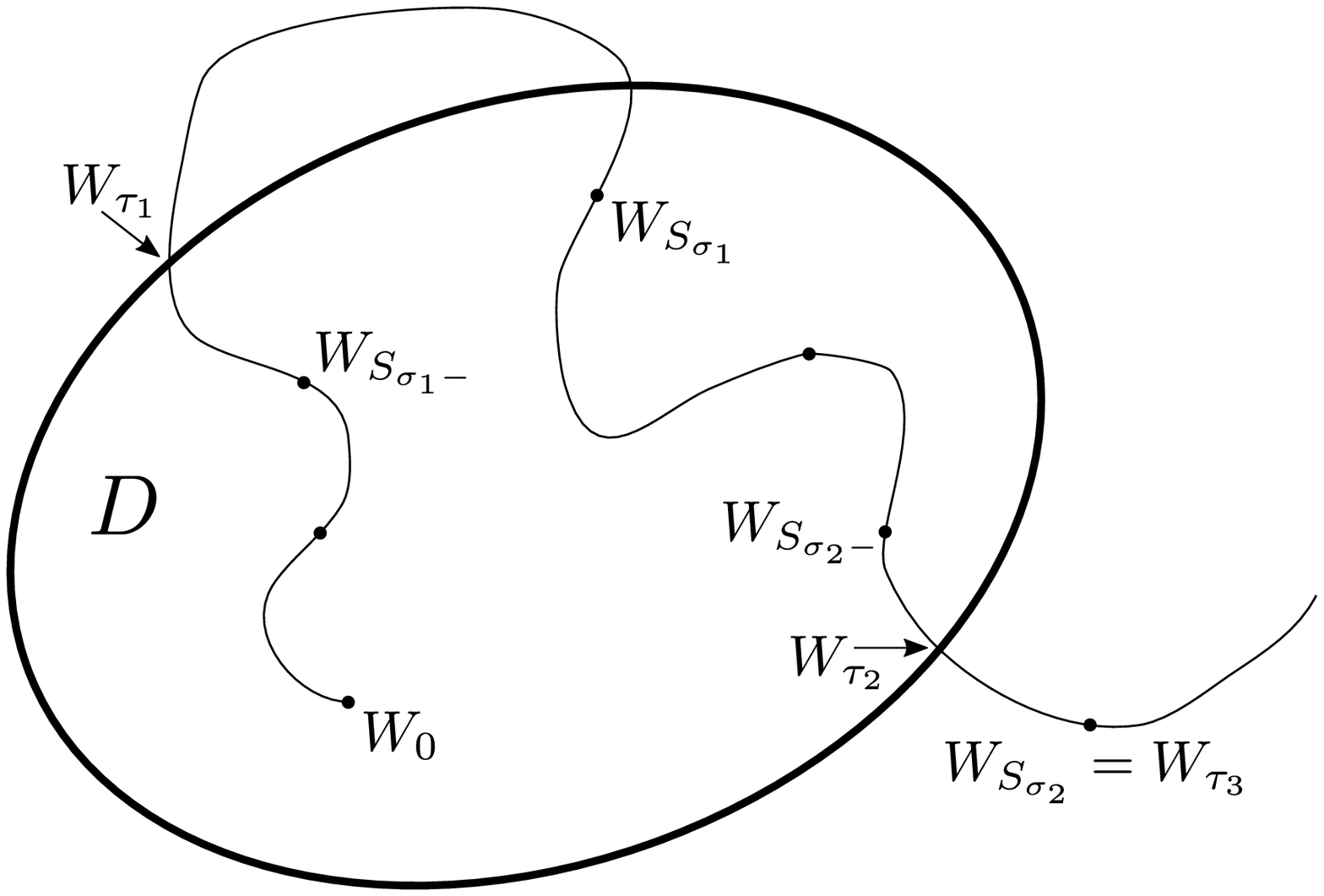}
\caption{A caricature of a Brownian path subordinated by a Poisson process where the dots indicate the range of the subordinate process. $N=2$ in this realization so exactly one resurrection of $Z^D$ is required to construct $\tau_D^Y$ as $\sigma_2$.}
\label{fig:subordinate}
\end{figure}

Let $\tau_0=0$ and $\sigma_0=0$ and define $\tau_{n+1}$ and $\sigma_{n+1}$ recursively by 
\[
\tau_{n+1}=\inf\{t\geq S_{\sigma_n}:W_t\notin D\}
\]
and
\[
\sigma_{n+1}=\inf\{t\geq 0:S_t\geq\tau_{n+1}\}.
\]
Notice that $\tau_1=\tau_D^W$ and $\sigma_1=\tau_D^Z$ and that the interlacing property $\tau_n\leq S_{\sigma_n}\leq \tau_{n+1}$ holds for all $n\geq 0$. Since the exit times of Brownian motion from a bounded domain as well as the first passage times of an unkilled subordinator across a level are all finite almost surely, it follows inductively that each of these stopping times is finite almost surely. Additionally, if for some $n\geq 0$ we have
\begin{equation}\label{eq:resurrect_end}
\tau_{n+1}=S_{\sigma_n},
\end{equation}
then $\sigma_{m+1}=\sigma_m$ and $\tau_{m+2}=\tau_{m+1}$ for all $m\geq n$. See Figure \ref{fig:subordinate} for an illustration when $S$ is a Poisson process.

Next we address the question of whether \eqref{eq:resurrect_end} holds for some $n$ almost surely. This is answered in the affirmative by showing that 
\begin{equation}\label{eq:N_def}
N:=\inf\{n\geq 0:\tau_{n+1}=S_{\sigma_n}\}
\end{equation}
is stochastically dominated by a geometric random variable. Towards this end, notice that for all $n$ we have
\begin{equation}\label{eq:stay_set}
\{\tau_{n+1}=S_{\sigma_n}\}=\{W_{S_{\sigma_n}}\notin D\}.
\end{equation}
The $\supset$ inclusion is immediate and the $\subset$ follows from the right-continuity of $W$ and the fact that $D^c$ is a closed set. 

Letting $A_n=\{W_{S_{\sigma_n}}\in D\}=\{Y_{\sigma_n}\in D\}$, we can use \eqref{eq:stay_set} along with the strong Markov property to write
\begin{align}
\mathbb{P}_x(\tau_{n+1}\neq S_{\sigma_n})&=\mathbb{P}_x\left(\tau_{2}\neq S_{\sigma_1}\cap\cdots\cap\tau_{n+1}\neq S_{\sigma_n}\right)\nonumber \\
&= \mathbb{E}_x\left[\prod_{j=1}^n \mathbbm{1}_{A_j}\right]\nonumber \\
&=\mathbb{E}_x\left[\prod_{j=1}^{n-1} \mathbbm{1}_{A_j}~\mathbb{P}_{Y_{\sigma_{n-1}}}\left(Y_{\sigma_1}\in D\right)\right]\nonumber \\
&\leq \mathbb{P}_x(\tau_{n}\neq S_{\sigma_{n-1}})\sup_{x\in \mathbb{R}^d}\mathbb{P}_x\left(W_{S_{\sigma_1}}\in D\right).\label{eq:stay_prob}
\end{align}

Since $\mathbb{P}_x\left(W_{S_{\sigma_1}}\in D\right)=0$ when $x\notin D$, we can restrict the supremum appearing in \eqref{eq:stay_prob} to $x\in D$ without affecting the inequality. In this case where we start $W$ from inside $D$, it follows from the continuity of $W$ that $W_{\tau_1}\in \partial D$. Assuming that $D$ is convex, there is a supporting hyperplane containing the point $W_{\tau_1}$ that divides $\mathbb{R}^d$ into two half-spaces, one which contains $D$ and another contained in the complement of $D$. Recalling that $\tau_1\leq S_{\sigma_1}$ by the interlacing property, now the rotational invariance of $W$ can be seen to imply that 
\begin{equation}\label{eq:stay_bound}
\sup_{x\in D}\mathbb{P}_x\left(W_{S_{\sigma_1}}\in D\right)\leq \frac{1}{2}.
\end{equation}
Combining \eqref{eq:stay_set}, \eqref{eq:stay_prob}, and \eqref{eq:stay_bound} with an inductive argument yields for all $x\in D$  
\begin{equation}\label{eq:stay_geo}
\left.\begin{aligned}
&\mathbb{P}_x\left(\tau_{n+1}\neq S_{\sigma_n}\right)\\
&\mathbb{P}_x\left(W_{S_{\sigma_n}}\in D\right)\\
&\mathbb{P}_x\left(Y_{\sigma_n}\in D\right)
\end{aligned}\right\}
\leq \frac{1}{2^n}.
\end{equation}

The estimate \eqref{eq:stay_geo} can be used with the definition of $N$ \eqref{eq:N_def} to conclude that
\begin{equation*}
\mathbb{P}_x(N>n)\leq \frac{1}{2^n}.
\end{equation*}
Hence for any starting point, it follows that $N$ is stochastically dominated by the geometric random variable $G$ with probability mass function
\begin{equation}\label{eq:geo}
\mathbb{P}(G=n)=\frac{1}{2^n},~n=1,2,\dots
\end{equation}

If instead of convexity we have an exterior cone condition on $D$, then \eqref{eq:stay_bound} will still hold but with the $\frac{1}{2}$ being replaced by some nonexplicit constant $\frac{1}{2}\leq C<1$; see the paragraph after the proof of Proposition 2.1 in \cite{killed_subordinate_domain} and the proof of Proposition 4.2 in \cite{two_sided} for similar considerations. Hence $N$ will still be stochastically dominated by a geometric random variable, in this case with mean $\frac{1}{1-C}$.

Next we relate $\sigma_N$ to $\tau_D^Y$. Since $Y_{\sigma_N}=W_{S_{\sigma_N}}\notin D$, we know that $\tau_D^Y\leq \sigma_N$. While this is sufficient for the upper bound we are most interested in, it also happens that the opposite inequality holds. To see that this is true, begin with the observation that as long as $W$ starts inside $D$, then for each $1\leq n\leq N$ we have $W_t\in D$ for all $t\in [S_{\sigma_{n-1}},S_{\sigma_n -})$. This implies $Y_t\in D$ for all $t\in [\sigma_{n-1},\sigma_n)$ for each $1\leq n\leq N$ hence $\tau_D^Y\geq \sigma_N$. If $W$ starts from outside $D$, then $N=0$ and $\tau_D^Y=0$ so $\tau_D^Y=\sigma_N$ trivially. Thus in either case we have
\begin{equation}\label{eq:sigma_rep}
\tau_D^Y=\sigma_N
\end{equation}
which can also be seen as a consequence of Proposition 3.2 in \cite{killed_subordinate}. It follows that $N-1$ corresponds to the number of resurrections required in order to construct $\tau_D^Y$ as $\sigma_N$ provided that $W$ starts in $D$; see Figure \ref{fig:subordinate} for an example where $N=2$.

\subsection{\texorpdfstring{Bounding $\|u_D^Y\|_\infty$ from above}{Bounding the Y torsion function from above}}

Finally, we use the above results to bound $\|u_D^Y\|_\infty$ in terms of the supremum of $\mathbb{E}_x\left[\tau_D^Z\right]$ taken over all starting points $x\in D$.
\begin{lem}\label{lem:resurrect}
Suppose $S$ is an unkilled subordinator and let $D\subset\mathbb{R}^d$ be a bounded convex domain. Then we have
\[
\left\|u_D^Y\right\|_\infty\leq 2\sup_{x\in D}\mathbb{E}_x\left[\tau_D^Z\right].
\]
\end{lem}

\begin{proof}
Let $A_n=\{Y_{\sigma_n}\in D\}$ as before and note that $\mathbbm{1}_{A_n}=\mathbbm{1}_{\{n<N\}}$ by \eqref{eq:N_def} and \eqref{eq:stay_set}. Hence starting from \eqref{eq:sigma_rep}, we can use the strong Markov property to write
\begin{equation*}
\mathbb{E}_x\left[\tau_D^Y\right]=\sum_{n=1}^\infty \mathbb{E}_x\left[(\sigma_n-\sigma_{n-1})\mathbbm{1}_{A_{n-1}}\right]=\sum_{n=1}^\infty \mathbb{E}_x\left[\mathbb{E}_{Y_{\sigma_{n-1}}}\left[\sigma_1\right]\mathbbm{1}_{A_{n-1}}\right]
\end{equation*}
for each $x\in D$. Now it follows from the estimate \eqref{eq:stay_geo} that
\begin{equation*}
\mathbb{E}_x\left[\tau_D^Y\right]\leq \sum_{n=1}^\infty \frac{1}{2^{n-1}}\sup_{y\in D}\mathbb{E}_y\left[\sigma_1\right]=2\sup_{y\in D}\mathbb{E}_y\left[\tau_D^Z\right].
\end{equation*}
\end{proof}

\section{Proofs of the main results}

\subsection{Proof of Theorem \ref{thm:torsion_bounds}}

\begin{proof}[Proof of Theorem \ref{thm:torsion_bounds}]
First we prove \eqref{eq:CS_torsion} under an additional boundedness assumption on $D$. Starting from $x\in D$, we know $W_t\in D$ for all $0\leq t< \tau_D^W$. Hence $Y_t=W_{S_t}\in D$ for all $0\leq t< \tau_D^Z$. This implies $\tau_D^Z\leq \tau_D^Y$, so under the hypotheses of Theorem \ref{thm:torsion_bounds} and the boundedness assumption, we can use Lemma \ref{lem:Jensen} to write
\begin{equation*}
\left\|u_D^Y\right\|_\infty \geq \sup_{x\in D}\mathbb{E}_x\left[\tau_D^Z\right]\geq\frac{1}{\Phi\left(\lambda_D^W\right)}.
\end{equation*}
Likewise, under the hypotheses of Theorem \ref{thm:torsion_bounds} and the boundedness assumption, we can apply both Lemmas \ref{lem:resurrect} and \ref{lem:Jensen} to conclude that
\begin{equation*}
\left\|u_D^Y\right\|_\infty\leq 2\sup_{x\in D}\mathbb{E}_x\left[\tau_D^Z\right]\leq 2\,V\Big(\left\|u_D^W\right\|_\infty\Big).
\end{equation*}
This proves Theorem \ref{thm:torsion_bounds} under the additional boundedness assumption.

Next we remove the boundedness assumption on $D$ by way of a localization argument borrowed from \cite{vdB_1,improved_Vogt}. Without loss of generality, assume that $0\in D$ and consider the sets $D_n:=D\cap B_n$ where $B_n$ is the ball of radius $n\in\mathbb{N}$ centered at the origin. Then $D_n$ is a bounded convex domain and $D_n\subset D_{n+1}$ for all $n$ with $\bigcup_{n\geq 1}=D$. Furthermore, it is clear from \eqref{eq:Rayleigh_quotient} that we have $\lambda_{D_1}^W\geq\lambda_{D_2}^W\geq\dots\geq\lambda_D^W$. Let $\{\varphi_m\}_{m\geq 1}$ be a sequence of functions in $C_c^\infty(D)$ with $\|\varphi_m\|_2=1$ such that $\int_D|\nabla \varphi_m|^2\D{x}\to\lambda_D^W$ as $m\to\infty$. Since $\supp\varphi_m\subset D_n$ for $n$ large enough, it follows that
\begin{equation}\label{eq:lambda_convergence}
\lim_{n\to\infty}\lambda_{D_n}^W\leq\int_D|\nabla \varphi_m|^2\D{x}.
\end{equation}
Now letting $m\to\infty$ in \eqref{eq:lambda_convergence} allows us to conclude that $\lambda_{D_n}^W\searrow\lambda_D^W$ as $n\to\infty$. 

We claim that $\|u_{D_n}^W\|_\infty\nearrow\|u_D^W\|_\infty$ as $n\to\infty$ also. Similarly to the principal eigenvalues, we have $\|u_{D_1}^W\|_\infty\leq\|u_{D_2}^W\|_\infty\leq\dots\leq\|u_D^W\|_\infty$. To see that convergence holds, let $\{x_m\}_{m\geq 1}$ be a sequence of points in $D$ such that $\mathbb{E}_{x_m}\left[\tau_D^W\right]\to\|u_D^W\|_\infty$ as $m\to\infty$. Since $x_m\in D_n$ for $n$ large enough, monotone convergence implies that
\begin{equation}\label{eq:torsion_convergence}
\lim_{n\to\infty}\left\|u_{D_n}^W\right\|_\infty\geq\lim_{n\to\infty}\mathbb{E}_{x_m}\left[\tau_{D_n}^W\right]=\mathbb{E}_{x_m}\left[\tau_D^W\right].
\end{equation}
Letting $m\to\infty$ in \eqref{eq:torsion_convergence} proves the claim. Furthermore, the same argument works for the subordinate Brownian motion so we also have $\|u_{D_n}^Y\|_\infty\nearrow\|u_D^Y\|_\infty$ as $n\to\infty$.

To finally prove Theorem \ref{thm:torsion_bounds}, we start with the localized version of \eqref{eq:CS_torsion}
\[
\frac{1}{\Phi\left(\lambda_{D_n}^W\right)}\leq\left\|u_{D_n}^Y\right\|_\infty\leq 2\,V\Big(\left\|u_{D_n}^W\right\|_\infty\Big),
\]
which holds since $D_n$ is a bounded convex domain, and then let $n\to\infty$ while noting that both $\Phi$ and $V$ are increasing and continuous on $[0,\infty)$ with $\Phi(0)=0$ and $\lim_{x\to\infty}V(x)=\infty$.
\end{proof}

\subsection{Proof of Theorem \ref{thm:main}}

Before proving the theorem, we first recall some important facts about stable subordinators. More specifically, we note that for $0<\alpha\leq 2$, an $\frac{\alpha}{2}$-stable subordinator has Laplace exponent $\Phi(\lambda)=\lambda^{\alpha/2}$ and is unkilled since $\mathbbm{k}=\Phi(0)=0$. Hence it follows from \eqref{eq:stable} and \eqref{eq:subordinate} that a $d$-dimensional Brownian motion subordinated by an $\frac{\alpha}{2}$-stable subordinator is a $d$-dimensional symmetric $\alpha$-stable process. In other words, if $S$ is an $\frac{\alpha}{2}$-stable subordinator, then the subordinate Brownian motion $Y$ is a $d$-dimensional symmetric $\alpha$-stable process. In particular, the fact that $\lambda_D^X=\lambda_D^Y$ and $\|u_D^X\|_\infty=\|u_D^Y\|_\infty$ is essential in the proof of Theorem \ref{thm:main} below. 

\begin{proof}[Proof of Theorem \ref{thm:main}]
We already know that \eqref{eq:main} holds for $\alpha=2$ by using $c_d$ and $C_d$ in \eqref{eq:spectral_bound}, so assume $0<\alpha<2$. Let $S$ be an $\frac{\alpha}{2}$-stable subordinator. We see from \cite[Section 5.2.2]{Bogdan_book} that $S$ has infinite L\'{e}vy measure and that its renewal function is 
\begin{equation}\label{eq:renewal_function}
V(x)=\frac{2}{\alpha\,\Gamma(\alpha/2)}x^{\alpha/2}.
\end{equation} 
Moreover, since $\Phi(\lambda)=\lambda^{\alpha/2}$, its conjugate Laplace exponent $\Phi^*(\lambda)=\lambda^{1-\alpha/2}$ is the Laplace exponent of a $\frac{2-\alpha}{2}$-stable subordinator. Hence $S$ is a special subordinator. 

Similarly to the proof of Theorem \ref{thm:torsion_bounds}, we first prove \eqref{eq:main} under an additional boundedness assumption on $D$. In this case the lower bound is simply a consequence of \eqref{eq:Smits}. We get the upper bound by using Theorem \ref{thm:torsion_bounds} with \eqref{eq:renewal_function}, then $C_d$ in \eqref{eq:spectral_bound}, and finally \eqref{eq:two_sided} to write
\begin{align*}
\left\|u_D^Y\right\|_\infty&\leq \frac{4}{\alpha\,\Gamma(\alpha/2)}\left\|u_D^W\right\|_\infty^{\alpha/2}\\
&\leq\frac{4}{\alpha\,\Gamma(\alpha/2)}\left(\frac{C_d}{\lambda_D^W}\right)^{\alpha/2}\\
&\leq\frac{4}{\alpha\,\Gamma(\alpha/2)}\frac{C_d^{\alpha/2}}{\lambda_D^Y}.
\end{align*}
Now the fact that the subordinate Brownian motion $Y$ is a $d$-dimensional symmetric $\alpha$-stable process proves Theorem \ref{thm:main} under the additional boundedness assumption.

We can remove the boundedness assumption with the same localization procedure used in the proof of Theorem \ref{thm:torsion_bounds}. The only difference is that we must now replace the Dirichlet form $\mathcal{E}_D^W(u,v):=\int_D\nabla u\cdot \nabla v\D{x}$ appearing in \eqref{eq:Rayleigh_quotient} with the appropriate quadratic form $\mathcal{E}_D^X$ associated to the Dirichlet fractional Laplacian. See \cite[Equation 84]{fractional_survey} for a description of $\mathcal{E}_D^X$ and also \cite[Section 5.1]{Frank} for more details. Now the same arguments show that $\lambda_{D_n}^X\searrow\lambda_D^X$ and $\|u_{D_n}^X\|_\infty\nearrow\|u_D^X\|_\infty$ as $n\to\infty$. Hence we can prove Theorem \ref{thm:main} by starting with the localized version of \eqref{eq:main} 
\[
\frac{1}{\lambda_{D_n}^X}\leq\left\|u_{D_n}^X\right\|_\infty\leq\frac{4}{\alpha\,\Gamma(\alpha/2)}\frac{C_d^{\alpha/2}}{\lambda_{D_n}^X},
\]
which holds since $D_n$ is a bounded convex domain, and then letting $n\to\infty$.
\end{proof}

\subsection{Proof of Lemma \ref{lem:comparison}}\label{sec:proofs}

\begin{proof}[Proof of Lemma \ref{lem:comparison}]
Tricomi's well-known asymptotic for Bessel zeros \cite{Tricomi} together with a lower bound of Lorch \cite{Lorch} imply that there exists $C\in (0,\infty)$ such that
\begin{equation}\label{eq:Bessel}
\frac{d^2}{4}\leq \lambda_B^W \leq \frac{d^2}{4}+C\,d^{4/3}
\end{equation}
holds for all $d=1,2,\dots$. Moreover, from \eqref{eq:two_sided}, we have
\begin{equation}\label{eq:two_sides}
\frac{1}{2}\left(\lambda_B^W\right)^{\alpha/2}\leq \lambda_B^X \leq \left(\lambda_B^W\right)^{\alpha/2}.
\end{equation}
Combining \eqref{eq:Bessel} with \eqref{eq:two_sides} and using the subadditivity of $x\mapsto x^{\alpha/2}$ leads to
\begin{equation}\label{eq:lambda_bound}
\frac{1}{2^{\alpha+1}}d^\alpha \leq \lambda_B^X \leq \frac{1}{2^\alpha}d^\alpha+C^{\alpha/2}d^{2\alpha/3}.
\end{equation}

The torsion function $u_B^X$ has the following expression \cite{ball_torsion}
\[
u_B^X(x)=\frac{\Gamma(d/2)\left(1-|x|^2\right)^{\alpha/2}}{2^\alpha \Gamma(1+\alpha/2)\Gamma(d/2+\alpha/2)},~|x|\leq 1
\]
from which we deduce that
\begin{equation}\label{eq:ball_sup}
\left\|u_B^X\right\|_\infty=\frac{\Gamma(d/2)}{2^\alpha \Gamma(1+\alpha/2)\Gamma(d/2+\alpha/2)}.
\end{equation}
Wendel's inequality for the ratio of gamma functions \cite[Equation 7]{Wendel} says that
\[
\left(\frac{x}{x+a}\right)^{1-a}\leq\frac{\Gamma(x+a)}{x^a\Gamma(x)}\leq 1 
\]
for any $x>0$ and $0<a<1$. Applied to the ratio appearing in \eqref{eq:ball_sup}, this gives
\[
\left(\frac{d}{2}\right)^{-\alpha/2}\leq\frac{\Gamma(d/2)}{\Gamma(d/2+\alpha/2)}\leq\left(\frac{d}{2}+\frac{\alpha}{2}\right)^{1-\alpha/2}\left(\frac{d}{2}\right)^{-1}.
\]
Applying this to \eqref{eq:ball_sup} and using the subadditivity of $x\mapsto x^{1-\alpha/2}$ results in
\begin{equation}\label{eq:tor_bound}
\frac{2^{-\alpha/2}}{\Gamma(1+\alpha/2)}d^{-\alpha/2}\leq\left\|u_B^X\right\|_\infty \leq \frac{2^{-\alpha/2}}{\Gamma(1+\alpha/2)}d^{-\alpha/2}+\frac{2^{-\alpha/2}\alpha^{1-\alpha/2}}{\Gamma(1+\alpha/2)}d^{-1}.
\end{equation}

The proof is completed by combining \eqref{eq:lambda_bound} with \eqref{eq:tor_bound} and absorbing the $d^{2\alpha/3-1}$ term that appears in the upper bound into the $d^{\alpha-1}$ or $d^{\alpha/6}$ term.
\end{proof}

\section*{Acknowledgments}

The author would like to thank Renming Song and Phanuel Mariano for their comments on an earlier draft. The author also acknowledges two anonymous referees whose many comments and suggestions led to a much improved paper.  


\end{document}